\documentclass[MSNbibl,number,citesort,seceqn,dvips]{arxbj}
\usepackage{upgreek}

\aid{0}
\volume{20}
\issue{3}
\pubyear{2014}
\firstpage{1620}
\lastpage{1646}
\doi{10.3150/13-BEJ535} 

\makeatletter
\newcommand{\RMo}{\mathrm{o}}
\newcommand{\RMO}{\mathrm{O}}

\newcommand{\RMe}{\mathrm{e}}

\newcommand{\RMi}{\mathrm{i}}

\newcommand{\mrmd}{\,\mathrm{d}}

\newtheorem{theorem}{Theorem}[section]
\newproclaim{definition}[theorem]{Definition}
\newtheorem{corollary}[theorem]{Corollary}
\newtheorem{proposition}[theorem]{Proposition}
\newremark{remark}[theorem]{Remark}
\newtheorem{lemma}[theorem]{Lemma}

\def\E{{\mathbf E}}
\def\R{{\mathbf R}}
\def\P{{\mathbf P}}

\def\Var{\operatorname{Var}}

\def\ep{\varepsilon}

\makeatother

\begin{document}
\begin{frontmatter}

\title{Fisher information and convergence to stable~laws}
\runtitle{Fisher information and stable laws}

\begin{aug}
\author[1]{\inits{S.G.}\fnms{S.G.} \snm{Bobkov}\thanksref{1}\ead[label=e1]{bobkov@math.umn.edu}},
\author[2]{\inits{G.P.}\fnms{G.P.} \snm{Chistyakov}\thanksref{2}\ead[label=e2,mark]{chistyak@math.uni-bielefeld.de}} \and
\author[2]{\inits{F.}\fnms{F.} \snm{G\"otze}\corref{}\thanksref{2}\ead[label=e3,mark]{goetze@mathematik.uni-bielefeld.de}}
\runauthor{S.G. Bobkov, G.P. Chistyakov and F. G\"otze} 
\address[1]{School of Mathematics, University of Minnesota,
127 Vincent Hall, 206 Church St. S.E., Minneapolis, MN 55455, USA.
\printead{e1}}
\address[2]{Faculty of Mathematics, University of Bielefeld,
Postfach 100131, 33501 Bielefeld, Germany.\\ \printead{e2,e3}}
\end{aug}

\received{\smonth{8} \syear{2012}}
\revised{\smonth{4} \syear{2013}}

%
\begin{abstract}
The convergence to stable laws is studied in relative Fisher information
for sums of i.i.d. random variables.
\end{abstract}

%

\end{frontmatter}

\section{Introduction}\label{sec1}

Let $(X_n)_{n \geq1}$ be independent identically distributed random
variables. Define the normalized sums
\[
Z_n = \frac{X_1 + \cdots+ X_n}{b_n} - a_n
\]
for given (non-random) normalizing sequences $a_n \in\R$ and $b_n >
0$. Assuming that $Z_n$ converges weakly in distribution to a random
variable $Z$ with a non-degenerate stable law, we consider the Fisher
information distance
\[
I(Z_n\|Z) = \int_{-\infty}^\infty\biggl(
\frac{p_n'(x)}{p_n(x)} - \frac{\psi'(x)}{\psi(x)} \biggr)^2 p_n(x)
\mrmd x,
\]
where $p_n$ and $\psi$ denote the densities of $Z_n$ and $Z$,
respectively. The definition makes sense, if $p_n$ is absolutely
continuous and is supported on the support interval of $\psi$, with
a Radon--Nikodym derivative $p_n'(x)$. Otherwise, put $
I(Z_n\|Z) = \infty$.

If $X_1$ has finite second moment with mean zero and variance one,
the classical central limit theorem is valid, that is,
$Z_n \Rightarrow Z$ weakly in distribution, with $a_n=0$,
$b_n = \sqrt{n}$, where $Z$ is standard normal.
In this case a striking result of Barron and Johnson \cite{B-J} indicates
that $I(Z_n\|Z) \rightarrow0$, as $n \rightarrow\infty$, as long as
$I(Z_n\|Z) < \infty$, for some $n$, that is, if for some $n$,
$Z_n$ has finite Fisher information
\[
I(Z_n) = \int_{-\infty}^\infty
\frac{p_n'(x)^2}{p_n(x)}\mrmd x.
\]

This observation considerably strengthens a number of results on the
central limit theorem for strong distances involving the total variation
and the relative entropy. It raises at the same time the question about
possible extensions to non-normal limit stable laws (as mentioned, e.g.,
in \cite{J}, page 104). The question turns out to be rather tricky, and
it is
not that evident that $I(Z_n)$ needs to be even bounded for large $n$
(a property which is guaranteed by Stam's inequality in case of
a finite second moment).

The present note gives an affirmative solution of the problem in case
of the so-called non-extremal stable laws, cf. Definition \ref{defi1.2} below.
In the sequel, we shall consider non-degenerate distributions, only.

%
\begin{theorem}\label{theo1.1}
Assume that the sequence of normalized sums
$Z_n$ defined above converges weakly in distribution to a random
variable $Z$ with a non-extremal stable limit law. Then
$I(Z_n\|Z) \rightarrow0$, as $n \rightarrow\infty$, if and only if
$I(Z_n\|Z) < \infty$ for some $n$.
\end{theorem}

The normal case is included in this assertion. Note, however, that if
$X_1$ has an infinite second moment, but still belongs to the domain
of normal attraction, we have $I(Z_n\|Z) = \infty$ for all $n$. Hence,
in this special case there is no convergence in relative Fisher
information.

In the remaining cases, $Z$ has a stable distribution with some parameters
$0 < \alpha< 2$, \mbox{$-1 \leq\beta\leq1$}, with characteristic function
$f(t) = \E\RMe^{\RMi tZ}$ described by
%
%
\begin{equation}\label{equ1.1}
f(t) = \exp\bigl\{\RMi at - c|t|^\alpha\bigl(1 + \RMi\beta\operatorname
{sign}(t)
\omega(t,\alpha) \bigr) \bigr\},
\end{equation}
where $a \in\R$, $c>0$, and $\omega(t,\alpha) = \tan(\frac{\uppi
\alpha}{2})$
in case $\alpha\neq1$, and $\omega(t,\alpha) = \frac{2}{\uppi}
\log|t|$
for $\alpha= 1$.
In particular, $|f(t)| = \RMe^{-c|t|^\alpha}$ which implies that
$Z$ has a smooth density $\psi(x)$.

%
\begin{definition}\label{defi1.2}
A stable distribution is called non-extremal,
if it is normal or, if $0<\alpha<2$ and $-1 < \beta< 1$ in (\ref{equ1.1}).
\end{definition}

In the latter case, the density $\psi$ of $Z$ is known to satisfy
asymptotic relations
%
%
\begin{equation}\label{equ1.2}
\psi(x) \sim c_0 |x|^{-(1+\alpha)} \qquad(x \rightarrow-\infty),
\qquad\psi(x) \sim c_1 x^{-(1+\alpha)} \qquad(x \rightarrow\infty)
\end{equation}
with some constants $c_0,c_1>0$. Since any stable distribution
is also unimodal (cf. \cite{Y}), $\psi$ has to be positive on the whole
real line, as follows from (\ref{equ1.2}).

The property that $X_1$ belongs to the domain of attraction of a stable
law of index $0 < \alpha< 2$ may be expressed explicitly in terms of
the distribution function $F_1(x) = \P\{X_1 \leq x\}$. Namely, we have
$Z_n \Rightarrow Z$ with some $b_n > 0$ and $a_n \in\R$, if and only
if
%
%
\begin{eqnarray}\label{equ1.3}
F_1(x) & = & \bigl(c_0 + \RMo(1)\bigr) |x|^{-\alpha}
B\bigl(|x|\bigr)\qquad(x \rightarrow-\infty),
\\
\label{equ1.4}
1-F_1(x) & = & \bigl(c_1 + \RMo(1)\bigr) x^{-\alpha}
B(x) \qquad(x \rightarrow\infty)
\end{eqnarray}
for some constants $c_0,c_1 \geq0$ that are not both zero, and where
$B(x)$ is a slowly varying function in the sense of Karamata.
This description reflects a certain behaviour of the characteristic
function $f_1(t) = \E\RMe^{\RMi tX_1}$ near the origin (cf. \cite{I-L,Z}).

In connection with Theorem \ref{theo1.1}, let us note that a similar assertion
has recently been proved in \cite{B-C-G1} for the relative entropy
\[
D(Z_n\|Z) = \int_{-\infty}^\infty
p_n(x) \log\frac{p_n(x)}{\psi
(x)}\mrmd x,
\]
called also the Kullback--Leibler distance form the distribution of
$Z_n$ to the distribution of $Z$.
It is shown that $D(Z_n\|Z) \rightarrow0$, if and only if
$Z_n \Rightarrow Z$ and $D(Z_n\|Z) < \infty$ for some $n$.
In the normal case this result is due to Barron \cite{B}, which in turn
goes back to the work by Linnik \cite{L}, initiating an
information-theoretic approach to the central limit theorem.

To compare with other strong types of convergence, in the normal case
it is known that, if $\E X_1 = \E Z$ and
$\Var(X_1) = \Var(Z) = \sigma^2$, then
%
%
\begin{equation}\label{equ1.5}
\frac{\sigma^2}{2} I(Z_n\|Z) \geq D(Z_n\|Z) \geq
\frac{1}{2} \|F_n - \Phi\|_{\mathrm{TV}}^2,
\end{equation}
where $\|F_n - \Phi\|_{\mathrm{TV}}$ is the distance in total variation
norm between the distributions of $Z_n$ and $Z$ (denoted here by $F_n$
and $\Phi$, resp.). The first relation in (\ref{equ1.5}), due to Stam
\cite{St}, may be viewed as an information theoretic variant of Gross'
logarithmic Sobolev inequality for the Gaussian measure.
The second one is a particular case of the Pinsker-type inequality
in which normality of $Z$ has no special role \cite{P,C,K}. Hence,
the convergence to the normal law in Fisher information distance is
a stronger property than in total variation and even than in relative
entropy. The question of how the Fisher information and entropic
distances are related to each other with respect to other stable laws
does not seem to have been addressed in the literature. Apparently it
is a question about the existence of certain weak logarithmic Sobolev
inequalities for probability distributions with heavy tails, and we
do not touch it here. However, it is natural to conjecture that the
situation is similar as in the normal case via a suitable analogue of
(\ref{equ1.5}).

Another obvious question concerns the description of distributions
satisfying the conditions of Theorem \ref{theo1.1}. In the non-normal case,
the property $I(Z_n\|Z) < \infty$ may be simplified to \mbox{$I(Z_n) <
\infty$}.
Taking, for example, $n = 1$, we obtain $I(X_1) < \infty$ as a
sufficient condition, which is however rather strong and may be considerably
weakened by choosing larger values of $n$. One may wonder therefore
what assumptions need to be added to (\ref{equ1.3})--(\ref{equ1.4}) in
terms of
$F_1$ or $f_1$ to obtain the convergence of $Z_n$ to $Z$ in relative Fisher
information. As shown in \cite{B-C-G3}, for some $n$, we have $I(Z_n)
< \infty$,
if and only if, for some $n$, $Z_n$ has a continuously differentiable
density $p_n$ such that
\[
\int_{-\infty}^\infty\bigl|p_n'(x)\bigr|
\mrmd x < \infty.
\]\eject\noindent
Still equivalently, for some $n$, $p_n$ has to be a function of bounded
variation. Moreover, if $X_1$ has a finite first absolute moment,
this property may be formulated explicitly in terms of the behaviour
of $f_1$ at infinity, as any of the following two equivalent assertions:

(a) For some $\ep>0$, $|f_1(t)| = \RMO(t^{-\ep})$, as
$t \rightarrow\infty$;

(b) For some $\nu> 0$,
%
%
\begin{equation}\label{equ1.6}
\int_{-\infty}^\infty\bigl|f_1(t)\bigr|^\nu
t^2 \mrmd t < \infty.
\end{equation}
This characterization may be used in Theorem \ref{theo1.1} in case
$1 < \alpha\leq2$, since then, by (\ref{equ1.3})--(\ref{equ1.4}), we have
$\E|X_1|^\delta< \infty$, for all $0 < \delta<\alpha$.

%
\begin{corollary}\label{corol1.3}
Assume that the sequence $Z_n$ as above
converges weakly in distribution to a random variable $Z$ with
a non-extremal stable limit law with index $1 < \alpha\leq2$.
Then $I(Z_n\|Z) \rightarrow0$, as $n \rightarrow\infty$, if and only
if (\ref{equ1.6}) holds for some $\nu> 0$.
\end{corollary}

In particular, this description is applicable to the usual central
limit theorem, that is, when $X_1$ has finite second moment. In this case
(cf. \cite{B-C-G3}), (\ref{equ1.6}) is equivalent to the formally
weaker condition
\[
\int_{-\infty}^\infty\bigl|f_1(t)\bigr|^\nu
|t| \mrmd t < \infty\qquad\mbox{for some } \nu> 0.
\]
However, removing the weight $|t|$ from the above integral, we obtain
an essentially weaker (so-called ``smoothness'') property
%
%
\begin{equation}\label{equ1.7}
\int_{-\infty}^\infty\bigl|f_1(t)\bigr|^\nu
\mrmd t < \infty\qquad\mbox{for some }\nu> 0.
\end{equation}
Once it is known that $Z_n \Rightarrow Z$ weakly in distribution
with a stable limit law (for the i.i.d. summands as above),
the condition (\ref{equ1.7}) allows one to strengthen the weak
convergence in the
following sense. It is equivalent to the property that, for some and
consequently for any sufficiently large~$n$, $Z_n$ has an absolutely
continuous distribution with a bounded continuous density $p_n$.
Moreover, in that and only that case, the uniform local limit theorem
holds: $\sup_x |p_n(x) - \psi(x)| \rightarrow0$,
as $n \rightarrow\infty$ (cf. \cite{I-L}).

The paper is organized as follows. First, we state some general bounds
on Fisher information and some properties of densities which can be
represented as convolutions of densities with finite Fisher information
(Sections~\ref{sec2}--\ref{sec4}). A main result used here has been
already proved in
recent work \cite{B-C-G3}. In Section~\ref{sec5}, we turn to the stable case
and discuss
a number of auxiliary results such as local limit theorems, as well as
questions about the behaviour of characteristic functions of $Z_n$ near
zero. In Section~\ref{sec6}, we reduce Theorem \ref{theo1.1} to showing
that the Fisher
information $I(Z_n)$ is bounded in $n$. The subsequent sections are
therefore focused on this boundedness problem. Section~\ref{sec7} introduces
a special decomposition of convolutions, and the final steps of
the proof of Theorem \ref{theo1.1} can be found in Section~\ref{sec8}.
We shall complement the proofs by comments explaining why the condition
(\ref{equ1.6}) is sufficient for the validity of Theorem \ref{theo1.1}.

\section{General results about Fisher information}\label{sec2}

%
\begin{definition}\label{defi2.1}
If a random variable $X$ has an absolutely
continuous density $p$ with Radon--Nikodym derivative $p'$, its Fisher
information is defined by
%
%
\begin{equation}\label{equ2.1}
I(X) = I(p) = \int_{\{p(x)>0\}} \frac{p'(x)^2}{p(x)} \mrmd x.
\end{equation}
In this case, if $\tilde p(x) = p(x)$ for almost all $x$ (a.e.),
put $I(\tilde p) = I(p)$. In any other case, \mbox{$I(X)=\infty$}.
\end{definition}

The equality (\ref{equ2.1}) appears as a particular case of the Fisher
information
\[
J(\theta) = \int_{-\infty}^\infty\biggl(
\frac{\partial p_\theta(x)}{\partial\theta} \biggr)^2 p_\theta(x)\mrmd x
\]
for the family of densities $p_\theta(x) = p(x-\theta)$
with respect to the location parameter $\theta\in\R$.

If $I(X)$ as defined in (\ref{equ2.1}) is finite, then necessarily the
distribution of $X$ has to be absolutely continuous with density $p(x)$
such that the derivative $p'(x)$ exists and is finite on a set of full
Lebesgue measure (and then $p$ will always be chosen to be a.e.
differentiable). Furthermore, one can show that, if $I(X) < \infty$,
then $p'(x) = 0$ at any point, where $p(x)=0$ (cf. \cite{B-C-G3}).
With this in mind, the integration in (\ref{equ2.1}) may be extended to the
whole real line.

It follows immediately from the definition that the $I$-functional is
translation invariant and homogeneous of order $-2$, that is,
$I(a+bX) = \frac{1}{b^2} I(X)$, for all $a \in\R$ and $b \neq0$.

Since the function $u^2/v$ is convex in the upper half-plane $u \in\R$,
$v>0$, this functional is convex. That is, for all densities
$p_1,\ldots,p_n$, we have Jensen's inequality
\[
I(\alpha_1 p_1 + \cdots+ \alpha_n
p_n) \leq\sum_{k=1}^n
\alpha_k I(p_k) \qquad\Biggl(\alpha_k > 0,
\sum_{k=1}^n \alpha_k = 1
\Biggr).
\]
The inequality may be generalized to arbitrary ``continuous'' mixtures
of densities. In particular, for the convolution
\[
p*q(x) = \int_{-\infty}^\infty p(x-y)q(y) \mrmd x
\]
of any two densities $p$ and $q$, we have
%
%
\begin{equation}\label{equ2.2}
I(p*q) \leq\min\bigl\{I(p),I(q)\bigr\}.
\end{equation}
In other words, if $X$ and $Y$ are independent random variables with these
densities, then
\[
I(X+Y) \leq\min\bigl\{I(X),I(Y)\bigr\}.
\]
This property may be viewed as monotonicity of the Fisher information:
this functional decreases when adding an independent summand.
In fact, a much stronger inequality is available.

%
\begin{proposition}[(Stam \cite{St})]\label{prop2.2}
If $X$ and $Y$ are independent
random variables, then
%
%
\begin{equation}\label{equ2.3}
\frac{1}{I(X+Y)} \geq\frac{1}{I(X)} + \frac{1}{I(Y)}.
\end{equation}
\end{proposition}

Let us also introduce the Fisher information distance
\[
I(X\|Z) = \int_{-\infty}^{+\infty} \biggl(\frac{p'(x)}{p(x)}
- \frac{\psi'(x)}{\psi(x)} \biggr)^2 p(x)\mrmd x
\]
with respect to a random variable $Z$ having a stable law.
We need the following elementary observation, which shows that the
question of boundedness of the Fisher information $I(Z_n)$ and of
the Fisher information distance $I(Z_n\|Z)$ for the normalized sums $Z_n$
as introduced in Theorem \ref{theo1.1} are in fact equivalent.

%
\begin{proposition}\label{prop2.3}
If $Z$ has a non-extremal stable law of
some index $0 < \alpha< 2$, then for any random variable $X$,
%
%
\begin{eqnarray}\label{equ2.4}
I(X\|Z) & \leq& 2 I(X) + c(Z),
\\
\label{equ2.5}
I(X) & \leq& 2 I(X\|Z) + c(Z),
\end{eqnarray}
where $c(Z)$ depends on the distribution of $Z$, only.
In particular, $I(X\|Z)<\infty$, if and only if $I(X)<\infty$.
\end{proposition}

\begin{pf}
The assertion is based on the fact that any non-extremal non-normal
stable distribution has a smooth positive density $\psi$ such that,
for all $k=1,2,\ldots\,$,
\[
\bigl|\bigl(\log\psi(x)\bigr)^{(k)} \bigr|\sim\frac{(k-1)!}{|x|^{k}} \qquad\bigl(|x|
\to
\infty\bigr)
\]
(cf. \cite{I-L,Z}). In particular,
$\frac{|\psi'(x)|}{\psi(x)} \sim\frac{1}{|x|}$, so
%
%
\begin{equation}\label{equ2.6}
\frac{|\psi'(x)|}{\psi(x)} \leq\frac{c}{1 + |x|} \qquad(x \in\R)
\end{equation}
with some positive constant $c$ (and the converse inequality is also
true with positive constant for all large $|x|$).
Hence, assuming that $I(X)<\infty$, then writing
\[
\biggl(\frac{p'(x)}{p(x)} - \frac{\psi'(x)}{\psi(x)} \biggr)^2 \leq2
\biggl(\frac{p'(x)}{p(x)} \biggr)^2 + 2 \biggl(\frac{\psi'(x)}{\psi
(x)}
\biggr)^2 \leq2 \biggl(\frac{p'(x)}{p(x)} \biggr)^2 +
2c^2
\]
and integrating this inequality with weight $p(x)$, we obtain (\ref{equ2.4}).
Similarly,
\[
\biggl(\frac{p'(x)}{p(x)} \biggr)^2 \leq2 \biggl(\frac{p'(x)}{p(x)}
- \frac{\psi'(x)}{\psi(x)} \biggr)^2 + 2c^2,
\]
which leads to (\ref{equ2.5}).
\end{pf}

Similar arguments for the normal case $(\alpha= 2)$ however lead
to a different conclusion. Indeed, if $Z \sim N(a,\sigma^2)$, we have
$\frac{\psi'(x)}{\psi(x)} = -\frac{x-a}{\sigma^2}$, and we get the
following proposition.

%
\begin{proposition}\label{prop2.4}
If $Z$ is normal, then
$I(X\|Z)<\infty$, if and only if $I(X) < \infty$ and $\E X^2 < \infty$.
\end{proposition}

Note that in case where $X$ and $Z$ have equal means and variances,
we have $I(X\|Z) = I(X) - I(Z)$.

\section{Connection with functions of bounded variation}\label{sec3}

Applying Cauchy's inequality and using the remark that
$p(x) = 0 \Rightarrow p'(x) = 0$ a.e., one immediately obtains from
Definition \ref{defi2.1} the following elementary lower bound on the
Fisher information.

%
\begin{proposition}\label{prop3.1}
If $X$ has an absolutely continuous density
$p$ with Radon--Nikodym derivative~$p'$, then
%
%
\begin{equation}\label{equ3.1}
\int_{-\infty}^\infty\bigl|p'(x)\bigr| \mrmd x \leq
\sqrt{I(X)}.
\end{equation}
\end{proposition}

Here, the integral represents the total variation norm of the function
$p$ as used in Real Analysis,
\[
\|p\|_{\mathrm{TV}} = \sup\sum_{k=1}^n
\bigl|p(x_k) - p(x_{k-1})\bigr|,
\]
where the supremum runs over all finite collections
$x_0 < x_1 < \cdots< x_n$.

The densities $p$ with finite total variation are vanishing at infinity
and are uniformly bounded by $\|p\|_{\mathrm{TV}}$. Moreover, their
characteristic functions
\[
f(t) = \int_{-\infty}^\infty\RMe^{\RMi tx} p(x) \mrmd x
\qquad(t \in\R)
\]
admit, by integration by parts, a simple upper bound
%
%
\begin{equation}\label{equ3.2}
\bigl|f(t)\bigr| \leq\frac{\|p\|_{\mathrm{TV}}}{|t|} \qquad(t \neq0).
\end{equation}
Hence, by Proposition \ref{prop3.1}, if a random variable $X$ has
finite Fisher
information, its density $p$ and characteristic function
$f(t) = \E\RMe^{\RMi tX}$ satisfy similar bounds
%
%
\begin{equation}\label{equ3.3}
\sup_x p(x) \leq\sqrt{I(X)}, \qquad\bigl|f(t)\bigr| \leq
\frac{\sqrt{I(X)}}{|t|} \qquad(t \neq0).
\end{equation}

In general, the inequality (\ref{equ3.1}) cannot be reversed, though
this is
possible for convolutions of three densities of bounded variation.
The following statement may be found in~\cite{B-C-G3}.

%
\begin{proposition}\label{prop3.2}
If independent random variables $X_j$
$(j=1,2,3)$ have densities $p_j$ of bounded variation, then
$S = X_1 + X_2 + X_3$ has finite Fisher information, and moreover,
%
%
\begin{equation}\label{equ3.4}
I(S) \leq\tfrac{1}{2} \bigl[ \|p_1\|_{\mathrm{TV}}
\|p_2\|_{\mathrm{TV}} + \|p_1\|_{\mathrm{TV}}
\|p_3\|_{\mathrm{TV}} + \|p_2\|_{\mathrm{TV}}
\|p_3\|_{\mathrm{TV}} \bigr].
\end{equation}
\end{proposition}

Note that the convolution of two densities of bounded variation
may have an infinite Fisher information. For example, the convolution
of the uniform distribution on $(-\frac{1}{2},\frac{1}{2})$ with itself
has the triangle density $p(x) = \max(1-|x|,0)$, in which case
$I(p) = \infty$.

%
\begin{remark}\label{rem3.3}
A similar bound on the Fisher information may also
be given in terms of characteristic functions. In view of (\ref
{equ3.4}), it
suffices to bound the total variation norm, and this can be done by
applying the inverse Fourier formula, at least in case of finite
first absolute moment. One can show that, if the characteristic function
$f(t)$ of a random variable $X$ is continuously differentiable
for $t>0$, and
%
%
\begin{equation}\label{equ3.5}
\int_{-\infty}^\infty t^2
\bigl(\bigl|f(t)\bigr|^2 + \bigl|f'(t)\bigr|^2 \bigr)\mrmd t <
\infty,
\end{equation}
then $X$ must have an absolutely continuous distribution with density
$p$ of bounded total variation satisfying
%
%
\begin{equation}\label{equ3.6}
\|p\|_{{\mathrm{TV}}} \leq\biggl(\int_{-\infty}^\infty\bigl|t
f(t)\bigr|^2\mrmd t \int_{-\infty}^\infty\bigl|\bigl(tf(t)
\bigr)'\bigr|^2\mrmd t \biggr)^{1/4}.
\end{equation}
We refer to \cite{B-C-G3} for details.
\end{remark}

\section{Classes of densities representable as convolutions}\label{sec4}

General bounds like (\ref{equ3.3}) may considerably be sharpened in the case
where $p$ is representable as convolution of several densities with
finite Fisher information. Here, we consider the collection
${\mathfrak P}_2(I)$ of all functions on the real line which can be
represented as convolution of two probability densities with Fisher
information at most $I$. Correspondingly, let
${\mathfrak P}_2 = \bigcup_I {\mathfrak P}_2(I)$ denote the collection of
all functions representable as convolution of two probability densities
with finite Fisher information. Note that, by (\ref{equ2.3}),
$I(p) \leq\frac{1}{2} I$, for any $p \in{\mathfrak P}_2(I)$.

Thus, a random variable $X = X_1 + X_2$ has density $p$ in
${\mathfrak P}_2$, if it may be written as
%
%
\begin{equation}\label{equ4.1}
p(x) = \int_{-\infty}^\infty p_1(x-y)
p_2(y)\mrmd x
\end{equation}
in terms of absolutely continuous densities $p_1,p_2$ of the
independent summands $X_1,X_2$ having finite Fisher information.
Differentiating under the integral sign,
we obtain a Radon--Nikodym derivative of the function $p$,
%
%
\begin{equation}\label{equ4.2}
p'(x) = \int_{-\infty}^\infty
p_1'(x-y) p_2(y)\mrmd y = \int
_{-\infty}^\infty p_1'(y)
p_2(x-y)\mrmd y.
\end{equation}
The latter expression shows that $p'$ is an absolutely continuous
function and has the Radon--Nikodym derivative
%
%
\begin{equation}\label{equ4.3}
p''(x) = \int_{-\infty}^\infty
p_1'(y) p_2'(x-y)\mrmd y.
\end{equation}
In other words, $p''$ appears as the convolution of the functions
$p_1'$ and $p_2'$ which are integrable, according to Proposition \ref{prop3.1}.

Note that equality (\ref{equ4.3}) defines $p''(x)$ at every individual
point $x$,
not just almost everywhere (which is typical for a Radon--Nikodym
derivative). Using the property $p_j(x) = 0 \Rightarrow p_j'(x) = 0$
in case of finite Fisher information, we obtain a similar implication
$p(x) = 0 \Rightarrow p''(x) = 0$, which holds for any $x$.

Moreover, since by (\ref{equ4.3}),
\[
\bigl|p''(x)\bigr| \leq\int_{-\infty}^\infty\bigl|p_1'(y)\bigr|
\bigl|p_2'(x-y)\bigr|\mrmd y,
\]
a direct application of the inequality (\ref{equ3.1}) together with Fubini's
theorem shows that $p'$ has finite total variation
\[
\bigl\|p'\bigr\|_{\mathrm{TV}} = \int_{-\infty}^\infty\bigl|p''(x)\bigr|
\mrmd x \leq I.
\]

These formulas may be used to derive various pointwise and integral
relations within the class ${\mathfrak P}_2$ such as the following
statement (which also summarizes the previous remarks).

%
\begin{proposition}\label{prop4.1}
Any density $p$ in ${\mathfrak P}_2(I)$
has an absolutely continuous derivative $p'$ of bounded variation
satisfying, for all $x \in\R$,
%
%
\begin{equation}\label{equ4.4}
\bigl|p'(x)\bigr| \leq I^{3/4} \sqrt{p(x)} \leq I.
\end{equation}
In addition,
%
%
\begin{equation}\label{equ4.5}
\int_{-\infty}^\infty\frac{p''(x)^2}{p(x)}\mrmd x \leq
I^2.
\end{equation}
\end{proposition}

To be more precise, integration in (\ref{equ4.5}) is restricted to the set
$\{p(x)>0\}$. This proposition can be found in \cite{B-C-G3}; since the
proof is short, we shall include it here for completeness.

\begin{pf*}{Proof of Proposition \ref{prop4.1}}
Starting with the representations (\ref{equ4.1})--(\ref{equ4.2}), in
which $I(p_j) \leq I$,
define the functions
$u_j(x) = \frac{p_j'(x)}{\sqrt{p_j(x)}}1_{\{p_j(x) > 0\}}$ ($j = 1,2$).
Applying Cauchy's inequality, we get
\begin{eqnarray*}
p'(x)^2 & = & \biggl(\int_{-\infty}^\infty
u_1(x-y) \cdot\sqrt{p_1(x-y)} p_2(y)\mrmd y
\biggr)^2
\\
& \leq& I(X_1) \int_{-\infty}^\infty
p_1(x-y) p_2(y)^2\mrmd y
\\
& \leq& I(X_1) \sqrt{I(X_2)} \int_{-\infty}^\infty
p_1(x-y) p_2(y)\mrmd y = I(X_1)
\sqrt{I(X_2)}p(x),
\end{eqnarray*}
where we used $p_2(y) \leq\sqrt{I(X_2)}$, according to (\ref{equ3.3}).
Hence, we obtain the first inequality in (\ref{equ4.4}), and the second follows
from $p(x) \leq\sqrt{I}$. Similarly, rewrite (\ref{equ4.3}) as
\[
p''(x) = \int_{-\infty}^\infty
\bigl(u_1(x-y) u_2(y) \bigr) \sqrt{p_1(x-y)p_2(y)}
\mrmd y
\]
to get
\[
p''(x)^2 \leq\int_{-\infty}^\infty
u_1(x-y)^2 u_2(y)^2\mrmd y \int
_{-\infty}^\infty p_1(x-y)p_2(y)
\mrmd y = u(x)^2 p(x),
\]
where we define $u \geq0$ by
\[
u(x)^2 = \int_{-\infty}^\infty
u_1(x-y)^2 u_2(y)^2\mrmd y.
\]
It follows that
\[
\int_{-\infty}^\infty u(x)^2\mrmd x =
I(X_1) I(X_2) \leq I^2,
\]
which implies (\ref{equ4.5}).
\end{pf*}

The analytic properties of densities in ${\mathfrak P}_2$ allow us
to make use of different formulas for the Fisher information
(by using integration by parts). For example,
\[
I(X) = -\int_{-\infty}^\infty p''(x)
\log p(x)\mrmd x,
\]
provided that the integrand is Lebesgue integrable.

We will need the following ``tail-type'' estimate for the Fisher
information.

%
\begin{corollary}\label{corol4.2}
If $p$ is in ${\mathfrak P}_2(I)$, then
for any $T$ real,
%
%
\begin{equation}\label{equ4.6}
\int_T^\infty\frac{p'(x)^2}{p(x)}\mrmd x \leq
I^{3/4} \sqrt{p(T)} \bigl|\log p(T)\bigr| + I \biggl(\int_T^\infty
p(x) \log^2 p(x)\mrmd x \biggr)^{1/2}.
\end{equation}
\end{corollary}

\begin{pf} Assuming that the last integral is finite, let us
decompose the open set $G = \{x>T\mbox{:}\allowbreak p(x) > 0\}$ into the union of
at most countably many disjoint intervals $(a_n,b_n)$,
$T \leq a_n < b_n \leq\infty$.

If $a_n>T$, we have $p(a_n) = 0$, so
$p'(x)\log p(x) \rightarrow0$, as $x \downarrow a_n$,
by Proposition \ref{prop4.1}. Similarly, $p(b_n) = 0$, if $b_n < \infty$,
and in addition $p(\infty) = 0$.

Let $a_n < T_1 < T_2 < b_n$. Since $p'$ is an absolutely continuous
function of bounded variation, integration by parts is justified and
yields
\[
\int_{T_1}^{T_2} \frac{p'(x)^2}{p(x)}\mrmd x = \int
_{T_1}^{T_2} p'(x)\mrmd\log p(x) =
p'(x) \log p(x) \bigg|_{x=T_1}^{T_2} - \int
_{T_1}^{T_2} p''(x) \log
p(x)\mrmd x.
\]
Letting $T_1 \rightarrow a_n$ and $T_2 \rightarrow b_n$, we get in case $a_n>T$
\[
\int_{a_n}^{b_n} \frac{p'(x)^2}{p(x)}\mrmd x = - \int
_{a_n}^{b_n} p''(x) \log
p(x)\mrmd x
\]
and
\[
\int_{a_n}^{b_n} \frac{p'(x)^2}{p(x)}\mrmd x =
-p'(T) \log p(T) - \int_{a_n}^{b_n}
p''(x) \log p(x)\mrmd x
\]
in case $a_n = T$ (if such $n$ exists). Anyhow, the summation over $n$
gives
%
%
\begin{equation}\label{equ4.7}
\int_G \frac{p'(x)^2}{p(x)}\mrmd x \leq\bigl|p'(T)
\log p(T)\bigr| + \int_G \bigl|p''(x)
\log p(x)\bigr|\mrmd x.
\end{equation}
Here the first term on the right-hand side can be estimated by virtue of
(\ref{equ4.4}), which leads to the first term on the right-hand side of
(\ref{equ4.6}).
Using (\ref{equ4.5}) together with Cauchy's inequality, for the last
integral we also have
\[
\biggl(\int_G \frac{|p''(x)|}{\sqrt{p(x)}} \sqrt{p(x)} \bigl|\log
p(x)\bigr|\mrmd x \biggr)^2 \leq I^2 \int_T^\infty
p(x) \log^2 p(x)\mrmd x,
\]
thus proving Corollary \ref{corol4.2}.
\end{pf}

\section{Stable laws and uniform local limit theorems}\label{sec5}

Let us return to the normalized sums
\[
Z_n = \frac{1}{b_n} (X_1 + \cdots+ X_n)
- a_n \qquad(a_n \in\R, b_n > 0),
\]
associated with independent identically distributed random variables
$(X_n)_{n \geq1}$. In this section, we discuss uniform limit theorems
for densities $p_n$ of $Z_n$ and behaviour of their characteristic
functions near the origin. As before, if $Z_n \Rightarrow Z$,
the density and the characteristic function of the stable limit $Z$
are denoted by $\psi$ and $f$, respectively.

Introduce the characteristic functions of $X_1$ and $Z_n$,
\[
f_1(t) = \E\RMe^{\RMi tX_1}, \qquad f_n(t) = \E
\RMe^{\RMi tZ_n} = \RMe^{-\RMi ta_n} f_1(t/b_n)^n
\qquad(t \in\R).
\]
To avoid confusion, we make the convention that $Z_1 = X_1$,
that is, $a_1 = 0$ and $b_1 = 1$.

%
\begin{proposition}\label{prop5.1}
Assume that $Z_n \Rightarrow Z$ weakly in
distribution. If
%
%
\begin{equation}\label{equ5.1}
\int_{-\infty}^\infty\bigl|f_1(t)\bigr|^\nu
\mrmd t < \infty\qquad\mbox{for some }\nu>0,
\end{equation}
then for all $n$ large enough, $Z_n$ have bounded continuous
densities $p_n$ such that
%
%
\begin{equation}\label{equ5.2}
\lim_{n\to\infty} \sup_x \bigl|p_n(x)-
\psi(x)\bigr| = 0.
\end{equation}
\end{proposition}

%
\begin{proposition}\label{prop5.2}
Assume that $Z_n \Rightarrow Z$ weakly in
distribution. If
%
%
\begin{equation}\label{equ5.3}
\int_{-\infty}^\infty\bigl|f_1(t)\bigr|^\nu
|t| \mrmd t < \infty \qquad\mbox{for some } \nu>0,
\end{equation}
then for all $n$ large enough, $Z_n$ have continuously
differentiable densities $p_n$ with bounded derivatives, and moreover
%
%
\begin{equation}\label{equ5.4}
\lim_{n\to\infty} \sup_x \bigl|p_n'(x)-
\psi'(x)\bigr| = 0.
\end{equation}
\end{proposition}

The first assertion is well known, cf. \cite{I-L}, page 126. The condition
(\ref{equ5.1}) is actually equivalent to the property that for all sufficiently
large $n$, say $n \geq n_0$, $Z_n$ have bounded continuous densities
$p_n$. In that case, the characteristic functions $f_n$ are integrable
whenever $n \geq2n_0$. Conversely, under (\ref{equ5.1}),
these densities for $n \geq\nu$ are given by the inversion formula
%
%
\begin{equation}\label{equ5.5}
p_n(x) = \frac{1}{2\uppi} \int_{-\infty}^\infty
\RMe^{-\RMi tx} f_n(t)\mrmd t.
\end{equation}

Under the stronger assumption (\ref{equ5.3}), the above equality may be
differentiated, and we get a similar representation for the derivative
%
%
\begin{equation}\label{equ5.6}
p_n'(x) = \frac{1}{2\uppi} \int_{-\infty}^\infty(-\RMi t)
\RMe^{-\RMi tx} f_n(t)\mrmd t.
\end{equation}

Although Proposition \ref{prop5.2} is not stated in \cite{I-L}, its
proof is similar
to the proof of Proposition~\ref{prop5.1}. An important ingredient in the
argument is the fact that the weak convergence $Z_n \Rightarrow Z$
forces $f_1$ to be regularly behaving near the origin. This fact can
also be used in the study of the boundedness of the Fisher information
distance $I(Z_n\|Z)$, so let us state it separately.

%
\begin{proposition}\label{prop5.3}
Let $Z_n \Rightarrow Z$ weakly in
distribution, where $Z$ has a stable law of index \mbox{$0 < \alpha< 2$}. Then
%
%
\begin{equation}\label{equ5.7}
\bigl|f_1(t)\bigr| = \exp\bigl\{-c |t|^{\alpha} h\bigl(1/|t|\bigr)\bigr\},
\end{equation}
where $c>0$ and $h(x)$ is a slowly varying function for
$x \rightarrow\infty$ such that
%
%
\begin{equation}\label{equ5.8}
\lim_{n \rightarrow\infty} \frac{nh(b_n)}{b_n^\alpha} = 1.
\end{equation}
Moreover, there is a constant $c>0$ such that, as $n \rightarrow\infty$,
%
%
\begin{equation}\label{equ5.9}
\P\bigl\{|X_1| > b_n\bigr\} \sim\frac{c}{n}.
\end{equation}
\end{proposition}

In comparison with (\ref{equ5.7}) a more precise statement is obtained
in \cite{I-L},
cf. Theorem 2.6.5, page~85. Namely, if $Z_n \Rightarrow Z$, where $Z$
has a stable distribution of index $0 < \alpha< 2$, then for all
$t$ small enough,
\[
f_1(t) = \exp\bigl\{\RMi\gamma t - c|t|^\alpha h\bigl(1/|t|\bigr)
\bigl(1 + \RMi\beta\operatorname{sign}(t) \omega(t,\alpha) \bigr) \bigr
\},
\]
where $\gamma$ is real, $c>0$, and the parameter $\beta\in[-1,1]$ and
the function $\omega(t,\alpha)$ are the same as in the representation
(\ref{equ1.1}) for the characteristic function $f(t)$ of $Z$. By lengthy
computations in the proof of Theorem 2.6.5 in \cite{I-L}, it was shown that
the function $B(x)$ appearing in the asymptotic relations (\ref
{equ1.3})--(\ref{equ1.4})
and the function $h(x)$ are connected via
\[
h(x) = \bigl(1 + \RMo(1)\bigr) B(x)\qquad\mbox{as }x \rightarrow\infty.
\]
Taking into account (\ref{equ5.8}), this yields (\ref{equ5.9}).

\begin{remark*}
As shown in \cite{I-L}, the representation (\ref{equ5.7})
together with
the relation (\ref{equ5.8}) remain to hold for $\alpha= 2$, that is,
when $Z$ is
normal. Note that, if $\E X_1^2 < \infty$, one may take $h(x)=1$ and
$b_n \sim\sqrt{n}$. In that case, $\P\{|X_1| > b_n\} = \RMo(\frac{1}{n})$,
as $n \rightarrow\infty$, so (\ref{equ5.9}) is no longer true.
\end{remark*}

Let us return to the local limit theorems.

\begin{pf*}{Proof of Proposition \ref{prop5.2}}
From (\ref{equ5.6}), we obtain the representation
\[
p_n'(x) - \psi'(x) = \frac{1}{2\uppi}
\int_{-\infty}^\infty(-\RMi t) \RMe^{-\RMi tx}
\bigl(f_n(t) - f(t)\bigr)\mrmd t.
\]
As is standard, we split the last integral into the three parts
$L_1$, $L_2$, $L_3$ corresponding to integration over the regions
$|t| \leq T_n$, $T_n < |t| < T_n'$ and $|t| \geq T_n'$, respectively.

By the weak convergence, $f_n(t) \rightarrow f(t)$ uniformly on
all intervals, and moreover,
\[
\delta_n = \max_{|t| \leq T_n} \bigl|f_n(t)- f(t)\bigr|
\rightarrow0\qquad\mbox{as } n \rightarrow\infty
\]
for some $T_n \rightarrow\infty$. Hence,
\[
|L_1| = \biggl|\int_{|t| \leq T_n} (-\RMi t) \RMe^{-\RMi tx}
\bigl(f_n(t) - f(t)\bigr)\mrmd t \biggr| \leq\delta_n
T_n^2 \rightarrow0,
\]
provided that $T_n$ grows to infinity sufficiently slowly (which may be
assumed).

Now, one of the consequences of (\ref{equ5.7}), using the above remark
about the
normal case, is that, given
$0 < \delta< \alpha$, the characteristic functions $f_n$ admit
on a relatively large interval the bound
%
%
\begin{equation}\label{equ5.10}
\bigl|f_n(t)\bigr| \leq\RMe^{-c(\delta) |t|^{\delta}} \qquad\bigl(|t| \leq\ep b_n\bigr)
\end{equation}
with some positive constants $\ep$ and $c(\delta)$ which are independent
of $n$, cf. \cite{I-L}, page 123. A similar bound holds for $f(t)$ itself,
which is also seen from the representation (\ref{equ1.1}).
Hence, choosing $T_n' = \ep b_n$, we have
\begin{eqnarray*}
|L_2| &=& \biggl|\int_{T_n < |t| < T_n'} (-\RMi t) \RMe^{-\RMi tx}
\bigl(f_n(t) - f(t)\bigr)\mrmd t \biggr| \\
&\leq&2\int_{|t| > T_n}
|t| \RMe^{-c(\delta)
|t|^{\delta}}\mrmd t \rightarrow0.
\end{eqnarray*}

Finally, put $c = \sup_{|t| \geq\ep} |f_1(t)|$. The condition (\ref{equ5.3})
ensures that $f_1(t) \rightarrow0$, as $t \rightarrow\infty$, so
$c<1$. Hence, for all $n \geq\nu$,
\begin{eqnarray*}
\int_{|t| \geq T_n'} |t| \bigl|f_n(t)\bigr|\mrmd t & = &
b_n^2 \int_{|t| \geq\ep} |t|
\bigl|f_1(t)\bigr|^n\mrmd t
\\
& \leq& b_n^2 c^{n-\nu} \int_{|t| \geq\ep}
|t| \bigl|f_1(t)\bigr|^\nu\mrmd t \rightarrow0.
\end{eqnarray*}
Thus, $L_3 \rightarrow0$, as well.
\end{pf*}

From (\ref{equ5.2}) and (\ref{equ5.4}), we immediately obtain the
convergence of
a ``truncated'' Fisher information distance.

%
\begin{corollary}\label{corol5.4}
Assume that $Z_n \Rightarrow Z$ weakly in
distribution, where $Z$ has a non-extremal stable law. If
$I(Z_{n_0})<\infty$ for some $n_0$, then for all $n$ large enough,
the random variables $Z_n$ admit continuously differentiable densities
$p_n$, and for every fixed $T>0$,
%
%
\begin{equation}\label{equ5.11}
\int_{-T}^T \biggl(\frac{p_n'(x)}{p_n(x)} -
\frac{\psi'(x)}{\psi(x)} \biggr)^2 p_n(x)\mrmd x = \RMo(1),\qquad n\to\infty.
\end{equation}
\end{corollary}

Recall that the densities $\psi$ of non-extremal stable laws are
everywhere positive, which is the only additional property needed
to show (\ref{equ5.11}) on the basis of (\ref{equ5.2}) and (\ref{equ5.4}).

Indeed, by the assumption, we have $I(Z_n)<\infty$, for all $n \geq n_0$,
and by (\ref{equ3.3}),
\[
\bigl|f_{n_0}(t)\bigr| \leq\frac{c}{|t|} \qquad(t \neq0)
\]
with $c = \sqrt{I(Z_{n_0})}$. Hence, the condition (\ref{equ5.3}) is fulfilled
with $\nu= 3n_0$. Therefore, we get both (\ref{equ5.2}) and (\ref
{equ5.4}), and
in particular, $p_n(x) \geq\ep> 0$ in $|x| \leq T$, for all $n$
large enough. As a result, the integrand in (\ref{equ5.11}) is uniformly
bounded from above by a sequence tending to zero.

\section{Moderate deviations}\label{sec6}

As before, for independent identically distributed random variables
$(X_n)_{n \geq1}$, put
%
%
\begin{equation}\label{equ6.1}
Z_n = \frac{X_1 + \cdots+ X_n}{b_n} - a_n \qquad(a_n \in
\R, b_n > 0).
\end{equation}

It is well known that if $Z_n \Rightarrow Z$, where $Z$ has a stable law
of some index $0 < \alpha\leq2$, then necessarily
%
%
\begin{equation}\label{equ6.2}
b_n = n^{1/\alpha}h(n),
\end{equation}
where $h$ is a slowly varying function in the sense of Karamata.

To study the behaviour of $I(Z_n\|Z)$ in the non-extremal non-normal
case, it is worthwhile noting that this Fisher information distance is
finite, if and only if $I(Z_n)$ is finite (Proposition~\ref{prop2.3}).
In the normal
case, $I(Z_n\|Z) < \infty$, if and only if $I(Z_n) < \infty$ and
$\E Z_n^2 < \infty$ (Proposition~\ref{prop2.4}). The latter is equivalent
to $\E X_1^2 < \infty$, and then for the weak convergence
$Z_n \Rightarrow Z$ with a standard normal limit one may take
$b_n = \sqrt{n \Var X_1 }$ and $a_n = \E X_1\sqrt n/\sqrt{\Var X_1}$.

In any case, the requirement that $I(Z_{n_0})<\infty$ implies that for
all $n \geq n_0$, $Z_n$ have absolutely continuous bounded densities
which we denote in the sequel by $p_n$. Moreover,
$p_n \in{\mathfrak P}_2$ whenever $n \geq2n_0$, and then, by
Proposition \ref{prop4.1}, $p_n$ have continuous derivatives $p_n'$ of bounded
variation.

As the next step towards Theorem \ref{theo1.1}, we prove the following lemma.

%
\begin{lemma}\label{lem6.1}
Assume that $Z_n \Rightarrow Z$ weakly in
distribution, where $Z$ has a non-extremal stable law. If
$\limsup_{n \rightarrow\infty} I(Z_n)<\infty$, then
%
%
\begin{equation}\label{equ6.3}
\lim_{n\to\infty} I(Z_n\|Z) = 0.
\end{equation}
\end{lemma}

\begin{pf} As before, denote by $\psi$ the density of $Z$, and put
$S_n = X_1 + \cdots+ X_n$.

By the assumptions, $I' = \sup_{n \geq n_0} I(Z_n) < \infty$
for some $n_0$, so
\[
I(S_n) \leq I' b_n^2, \qquad n
\geq n_0.
\]
If $n \geq2n_0$, write $n = n_1 + n_2$ with $n_1 = [\frac{n}{2}]$,
$n_2 = n - n_1$. Then $n_1 \geq n_0$ and $n_2 \geq n_0$, and hence
\[
I(S_{n_1}) \leq I' b_{n_1}^2 \leq
I b_n^2, \qquad I(S_n - S_{n_1})
\leq I' b_{n_2}^2 \leq I b_n^2
\]
with some constant $I$ in view of the almost polynomial behaviour of
$b_n$ as described in (\ref{equ6.2}). Thus,
\[
Z_n = \biggl(\frac{S_{n_1}}{b_n} - a_n \biggr) +
\frac{S_n - S_{n_1}}{b_n}
\]
represents the sum of two independent random variables with Fisher
information at most $I$. Therefore, $p_n \in{\mathfrak P}_2(I)$,
for all $n \geq2n_0$, and we may invoke Corollary \ref{corol4.2}.

In view of Corollary \ref{corol5.4} we only need to show that, given
$\ep> 0$,
one may choose $T>0$ such that the integral
\[
J = \int_{|x|>T} \biggl(\frac{p_n'(x)}{p_n(x)} - \frac{\psi'(x)}{\psi
(x)}
\biggr)^2p_n(x)\mrmd x
\]
is smaller than $\ep$, for all $n$ large enough.

Clearly, $J \leq2J_1 + 2J_2$, where
\[
J_1 = \int_{|x|>T}\frac{p_n'(x)^2}{p_n(x)}\mrmd x, \qquad
J_2 = \int_{|x|>T} \biggl(\frac{\psi'(x)}{\psi(x)}
\biggr)^2 p_n(x)\mrmd x.
\]

Recall that in case $0 < \alpha< 2$, we have
$\frac{|\psi'(x)|}{\psi(x)} \leq\frac{c}{1+|x|}$ with a constant $c$
depending on $\psi$, only (cf. (\ref{equ2.6})). Hence,
\[
J_2 \le\biggl(\frac{c}{1+T} \biggr)^2,
\]
which thus can be made as small, as we wish.

If $\alpha= 2$ and $\E X_1^2 < \infty$, assume without loss of
generality that $\E X_1 = 0$, $\E X_1^2 = 1$, so that $\psi$ is
a standard normal density, and
\[
J_2 = \int_{|x|>T} x^2
p_n(x)\mrmd x.
\]
To bound these integrals, we appeal to the well-known large deviation
relation
\[
\P\bigl\{|\xi| \geq T\bigr\} \leq T \int_0^{2/T} \bigl(1 -
\operatorname{Re} f(t)\bigr)\mrmd t,
\]
holding true for any random variable $\xi$ with characteristic function
$f(t)$. If $\E\xi^2 = 1$, and $F$ is the distribution function of
$\xi$,
one may apply the same bound to the probability measure $x^2 \mrmd
F(x)$ on
the real line, and then it yields
\[
\int_{|x| \geq T} x^2\mrmd F(x) \leq T \int
_0^{2/T} \bigl(1 + \operatorname{Re} f''(t)
\bigr)\mrmd t.
\]
Hence,
\[
J_2 \leq T \int_0^{2/T} \bigl(1 +
\operatorname{Re} f_n''(t)\bigr)\mrmd t,
\]
where $f_n$ denote the characteristic functions of $Z_n$. But,
letting $g(t) = \RMe^{-t^2/2}$, as a variant of the central limit theorem,
for any $c>0$, one has
$\sup_{|t| \leq c} |f_n''(t) - g''(t)| \rightarrow0$, as
$n \rightarrow\infty$, while $1 + g''(t) \rightarrow0$, as
$t \rightarrow0$. This shows that, for $T$ and $n$ large enough,
$J_2$ will be smaller than any prescribed positive number.

It remains to estimate $J_1$. We now apply (\ref{equ4.6}) giving
%
%
\begin{eqnarray}\label{equ6.4}
J_1 & \leq& I^{3/4} \bigl(\sqrt{p_n(T)} \bigl|\log
p_n(T)\bigr| + \sqrt{p_n(-T)} \bigl|\log p_n(-T)\bigr|
\bigr)
\nonumber\\[-8pt]\\[-8pt]
& &{} + 2 I \biggl(\int_{|x| \geq T} p_n(x)
\log^2 p_n(x)\mrmd x \biggr)^{1/2}.\nonumber
\end{eqnarray}
Using the uniform local limit theorem in the form (\ref{equ5.2})
together with
the asymptotic relation (\ref{equ1.2}) for $\psi(x)$ at infinity, we
easily get
%
%
\begin{equation}\label{equ6.5}
\sqrt{p_n(\pm T)} \bigl|\log p_n(\pm T)\bigr| \le c
\frac{\log T}{\sqrt{T}} + \ep_n,
\end{equation}
which holds for all sufficiently large $n$ and all $T \geq T_0$
with $\ep_n \rightarrow0$ (as $n \rightarrow\infty$) and with
constants $c > 0$ and $T_0$ large enough (depending on $\psi$, only).

To bound the integral in (\ref{equ6.4}), we partition $\{x\dvt |x|
\geq T\}$ into
the set
\[
A = \bigl\{x\dvt|x| \geq T, p_n(x) \le|x|^{-4}\bigr\}
\]
and its complement $B$. By the definition,
%
%
\begin{equation}\label{equ6.6}
\int_A p_n(x)\log^2
p_n(x)\mrmd x \le16 \int_{|x| \geq T} |x|^{-4}
\log^2 |x|\mrmd x \le\frac{32}{T}.
\end{equation}
On the other hand, $p_n$ are uniformly bounded, namely,
$\sup p_n(x)\le\sqrt{I}$, for all $n \geq2n_0$ (cf. (\ref{equ3.3})).
Hence, on the set
$B$,
\[
\bigl|\log p_n(x)\bigr| \leq\log\frac{\sqrt{I}}{p_n(x)} + |{\log\sqrt{I}}| \leq4
\log|x| + |{\log I}|
\]
and therefore
%
%
\begin{equation}\label{equ6.7}
\int_B p_n(x)\log^2
p_n(x)\mrmd x \le c\int_{|x| \geq T} p_n(x)
\log^2 |x|\mrmd x,
\end{equation}
where the constant depends on $I$.

Finally, we use the property that the moments $\E|Z_n|^\delta$ are
uniformly bounded in $n$, whenever $0 < \delta<\alpha$
(cf. \cite{I-L}, page 142). Choosing $\delta= \alpha/2$ and using
an elementary bound $|x|^{\alpha/4} \geq c_\alpha\log^2 |x|$
for $|x| \geq T_0$, we obtain with some constant $K$ that
\begin{eqnarray*}
K &\geq&\E|Z_n|^{\alpha/2} \geq T^{\alpha/4} \E
|Z_n|^{\alpha/4} 1_{\{|Z_n| \geq T\}}
\\
& = & T^{\alpha/4} \int_{|x| \geq T} |x|^{\alpha/4}
p_n(x)\mrmd x \\
&\geq& c_\alpha T^{\alpha/4} \int
_{|x| \geq T} p_n(x)\log^2 |x|\mrmd x.
\end{eqnarray*}

Thus, the second integral in (\ref{equ6.7}) may be bounded by
$cT^{-\alpha/4}$
with some constant $c$ independent of $n$. Combining this with (\ref{equ6.6}),
we obtain a similar bound for the integral in (\ref{equ6.4}), and
taking into
account (\ref{equ6.5}), we get $J_1 \leq cT^{-\alpha/8} + \ep_n$.
This completes the proof of Lemma \ref{lem6.1}.
\end{pf}

\section{Binomial decomposition of convolutions}\label{sec7}

To show that the assumption
$\limsup_{n \rightarrow\infty} I(Z_n)<\infty$ in Lemma \ref{lem6.1} holds
as long as $I(Z_{n_0}) < \infty$ for some $n_0$, we introduce a special
decomposition of densities of $Z_n$. It is needed for the case
$0<\alpha<2$, so this will be assumed below. Moreover, let
$Z_n \Rightarrow Z$ weakly in distribution, where $Z$ has a non-extremal
stable law with index $\alpha$.

To simplify the argument, assume $n_0 = 1$, so that
$I(p) = I(X_1) < \infty$, where $p$ denotes the density of $X_1$.
In fact, we only consider the shifted normalized sums
\[
\tilde Z_n = Z_n + a_n =
\frac{X_1 + \cdots+ X_n}{b_n},
\]
and for the notational convenience, denote their densities by $p_n$.
Note that, by the translation invariance, $I(Z_n) = I(\tilde Z_n)$.

Keeping the same notations as in the previous sections, we use a
suitable truncation (which is actually not needed in case $\alpha> 1$).
Introduce the probability densities
\[
\tilde{p}_n(x) = \frac{b_n}{1-\delta_n} p(b_n x)
1_{\{|x| \leq1\}}, \qquad\tilde{q}_n(x) = \frac{b_n}{\delta_n}
p(b_n x) 1_{\{|x| > 1\}}
\]
together with their characteristic functions
\[
\tilde{f}_n(t) = \frac{1}{1-\delta_n} \int_{-b_n}^{b_n}
\RMe^{\RMi tx/b_n} p(x)\mrmd x,\qquad \tilde{g}_n(t) = \frac{1}{\delta_n} \int
_{|x|>b_n} \RMe^{\RMi tx/b_n} p(x)\mrmd x,
\]
where $\delta_n = \int_{|x|>b_n}p(x)\mrmd x$. Recall that
$\delta_n \sim\frac{c}{n}$ with some constant $c>0$, as emphasized in
Proposition \ref{prop5.3}, cf. (\ref{equ5.9}).

Then we have a binomial decomposition for convolutions
%
%
\begin{equation}\label{equ7.1}
p_n = \bigl((1-\delta_n) \tilde{p}_n +
\delta_n\tilde{q}_n\bigr)^{n*} = \sum
_{k=0}^n \pmatrix{n \cr k} (1-\delta_n)^k
\delta_n^{n-k} \tilde{p}_n^{k*}*
\tilde{q}_n^{(n-k)*}.
\end{equation}
Note that each convolution $\tilde{p}_n^{k*}* \tilde{q}_n^{(n-k)*}$
appearing in this weighted sum represents a probability density with
characteristic function $\tilde f_n(t)^k \tilde g_n(t)^{n-k}$.

In this section, we establish some properties of $\tilde{f}_n$, which
will be needed in the proof of Theorem \ref{theo1.1}. The corresponding density
$\tilde p_n$ is supported on $[-1,1]$, however, it does not need
to have mean zero. So, put
\[
d_n = \int_{-1}^1 x
\tilde{p}_n(x)\mrmd x = \frac{1}{b_n (1-\delta_n)} \int_{-b_n}^{b_n}
x p(x)\mrmd x
\]
and define
\[
\psi_n(t) = \RMe^{-\RMi td_n} \tilde{f}_n(t),
\]
which is the characteristic function of the centered random variable
$\xi- d_n$, when $\xi$ has density~$\tilde{p}_n$. Thus,
$\psi_n$ corresponds to the density $r_n(x) = \tilde{p}_n(x+d_n)$,
with $\psi_n'(0)=0$.

The next two lemmas do not use the assumption $I(p) < \infty$ and
may be stated for general distributions from the domain of attraction
of these stable laws.

%
\begin{lemma}\label{lem7.1}
For all real $t$, with some constant $C$ depending
only on $p$,
%
%
\begin{equation}\label{equ7.2}
\bigl|\psi_n'(t)\bigr| \leq\frac{C}{n} |t|.
\end{equation}
\end{lemma}

\begin{pf} The characteristic function $\psi_n$ corresponds to the
density $\tilde p_n(x + d_n)$. Using the property $\psi_n'(0)=0$,
one may write
\begin{eqnarray*}
\psi_n'(t) & = & \int_{-1}^1
\RMi (x-d_n) \bigl(\RMe^{\RMi t(x-d_n)} - 1 \bigr) \tilde p_n(x)
\mrmd x
\\
& = & \frac{\RMi b_n}{1-\delta_n} \int_{-1}^1
(x-d_n) \bigl(\RMe^{\RMi t(x-d_n)} - 1 \bigr) p(b_n x)\mrmd x
\\
& = & \frac{i}{1-\delta_n} \int_{-b_n}^{b_n} \biggl(
\frac{x}{b_n}-d_n \biggr) \bigl(\RMe^{\RMi t({x}/{b_n} - d_n)} - 1
\bigr)
\mrmd F_1(x),
\end{eqnarray*}
where $F_1$ is the distribution function of $X_1$. Using
$|\RMe^{is} - 1| \leq|s|$ ($s \in\R$), we deduce obvious upper bounds
\begin{eqnarray*}
\bigl|\psi_n'(t)\bigr| & \le& \frac{|t|}{1-\delta_n} \int
_{-b_n}^{b_n} \biggl(\frac{x}{b_n}-d_n
\biggr)^2\mrmd F_1(x)
\\
& \le& \frac{2|t|}{b_n^2 (1-\delta_n)} \int_{-b_n}^{b_n}
x^2\mrmd F_1(x) + \frac{2|t|}{1-\delta_n} d_n^2.
\end{eqnarray*}
Integrating by parts, we have
\begin{eqnarray*}
\int_{-b_n}^{b_n} x^2
\mrmd F_1(x) & = & -b_n^2 \bigl(1-F_1(b_n)+F_1(-b_n)
\bigr) + 2\int_0^{b_n} x\bigl(1-F_1(x)+F_1(-x)
\bigr)\mrmd x
\\
& \le& 2\int_0^{b_n} x\bigl(1-F_1(x)+F_1(-x)
\bigr)\mrmd x
\end{eqnarray*}
and similarly
\begin{eqnarray*}
\int_{-b_n}^{b_n} |x|\mrmd F_1(x) & = &
-b_n \bigl(1-F_1(b_n)+F_1(-b_n)
\bigr) + \int_0^{b_n} \bigl(1-F_1(x)+F_1(-x)
\bigr)\mrmd x
\\
& \le& \int_0^{b_n} \bigl(1-F_1(x)+F_1(-x)
\bigr)\mrmd x.
\end{eqnarray*}
Since $1-\delta_n \rightarrow1$, we get
%
%
\begin{eqnarray}\label{equ7.3}
\bigl|\psi_n'(t)\bigr| & \leq& \frac{C|t|}{b_n^2} \int
_0^{b_n} x\bigl(1-F_1(x)+F_1(-x)
\bigr)\mrmd x
\nonumber\\[-8pt]\\[-8pt]
& &{} + \frac{C|t|}{b_n^2} \biggl(\int_0^{b_n}
\bigl(1-F_1(x)+F_1(-x)\bigr)\mrmd x \biggr)^2\nonumber
\end{eqnarray}
with some constant $C$ depending on $p$.

Recall that in the asymptotical formulas (\ref{equ1.3})--(\ref{equ1.4})
for $F_1$,
the function $B$ is equivalent to the slowly varying function $h$
associated with the characteristic function of $X_1$. Thus, with some
$c_0 \geq0$, $c_1 \geq0$ $(c_0 + c_1 > 0)$, we have
\[
F_1(x) = \frac{c_0+\RMo(1)}{(-x)^{\alpha}} h(-x),\qquad x<0;\qquad
F_1(x) = 1-\frac{c_1+\RMo(1)}{x^{\alpha}} h(x),\qquad x>0.
\]
Hence, up to a constant, the first integral in (\ref{equ7.3}) does not exceed
\[
\int_0^{b_n} \frac{h(x)}{x^{\alpha-1}}\mrmd x =
b_n^{2-\alpha} h(b_n)\int_0^1
\frac{h(sb_n)}{h(b_n)} \,\frac
{\mathrm{d}s}{s^{\alpha-1}}.
\]
But, by the well-known result on slowly varying functions
(\cite{Se}, pages 66--67),
\[
\int_0^1\frac{h(sb_n)}{h(b_n)} \,\frac{\mathrm{d}s}{s^{\alpha-1}}
\to\int_0^1\frac{\mathrm{d}s}{s^{\alpha-1}} =
\frac{1}{2-\alpha}\qquad\mbox{as } n\to\infty.
\]
Therefore, with some constants $C_1$, $C_2$,
\[
\frac{1}{b_n^2} \int_0^{b_n} x
\bigl(1-F_1(x)+F_1(-x)\bigr)\mrmd x \leq C_1
b_n^{-\alpha} h(b_n) \leq\frac{C_2}{n},
\]
where we have applied equation (\ref{equ5.8}) of Proposition \ref
{prop5.3}, telling us
that $h(b_n) \sim b_n^\alpha/n$.

Now, consider the second integral in (\ref{equ7.3}). In case $\alpha< 1$,
again by \cite{Se}, applied to the value $\alpha+1$,
\[
\int_0^1\frac{h(sb_n)}{h(b_n)} \,\frac{\mathrm{d}s}{s^\alpha}
\longrightarrow\int_0^1\frac{\mathrm{d}s}{s^\alpha} =
\frac{1}{1-\alpha}\qquad\mbox{as } n\to\infty.
\]
Hence, using the asymptotic for $F_1$, the second integral in (\ref{equ7.3})
does not exceed, up to a constant,
\[
\int_0^{b_n} \frac{h(x)}{x^\alpha}\mrmd x =
b_n^{1-\alpha} h(b_n) \int_0^1
\frac{h(sb_n)}{h(b_n)} \,\frac{\mathrm{d}s}{s^\alpha} \sim\frac{b_n}{(1-\alpha
) n}.
\]
As a result,
\[
\frac{1}{b_n^2} \biggl(\int_0^{b_n}
\bigl(1-F_1(x)+F_1(-x)\bigr)\mrmd x \biggr)^2
\leq\frac{C_3}{n^2}
\]
with some constant $C_3$, depending on $p$ and $\alpha$.

The case $1 < \alpha< 2$ is simpler, since then
\[
\int_0^\infty\bigl(1-F_1(x)+F_1(-x)
\bigr)\mrmd x < \infty,
\]
while the factor $\frac{1}{b_n^2}$ behaves like $n^{-2/\alpha}$
(up to a slowly growing sequence), so it decays faster than $1/n$.

Finally, in case $\alpha= 1$, using the bound $h(x) \leq C_\ep x^{\ep}$,
$x \geq1$ (where $\ep> 0$ is any prescribed number), we see that,
for large $n$ the second integral in (\ref{equ7.3}) does not exceed,
up to a constant,
\[
1 + \int_1^{b_n} \frac{h(x)}{x}\mrmd x \leq C
b_n^{1/4}.
\]
This yields
\[
\frac{1}{b_n^2}\biggl(\int_0^{b_n}
\bigl(1-F_1(x)+F_1(-x)\bigr)\mrmd x \biggr)^2
\leq\frac{C}{b_n^{3/2}}
\]
with some constant $C$ depending on the density $p$. But the ratio
$\frac{C}{b_n^{3/2}}$ behaves like $n^{-3/2}$ up to a slowly
growing sequence, so it decays faster than $\frac{1}{n}$, as well. Thus,
in all cases
\[
\frac{1}{b_n^2} \biggl(\int_0^{b_n}
\bigl(1-F_1(x)+F_1(-x)\bigr)\mrmd x \biggr)^2 = \RMO
\biggl(\frac{1}{n} \biggr).
\]
Lemma \ref{lem7.1} is proved.
\end{pf}

%
\begin{lemma}\label{lem7.2}
Let $\delta\in(0,\alpha)$ and $\eta\in(0,1)$
be fixed. There exist positive constants $\ep$, $c$, $C$, depending on
$p,\delta,\eta$, with the following property: if $k \geq\eta n$, then
%
%
\begin{equation}\label{equ7.4}
\bigl|\psi_n(t)\bigr|^k = \bigl|\tilde{f}_n(t)\bigr|^k
\leq C \RMe^{-c|t|^{\delta}}\qquad\mbox{for }|t| \leq\ep b_n.
\end{equation}
\end{lemma}

\begin{pf} This is an analogue of the bound (\ref{equ5.10}) for the
characteristic functions of $Z_n$. In order to prove this upper bound,
assume $|t| \geq1$ and note that
%
%
\begin{equation}\label{equ7.5}
\tilde{f}_n(t) = \frac{1}{1-\delta_n} \bigl(f_1(t/b_n)-
\delta_n \tilde{g}_n(t) \bigr), \qquad t \in\R.
\end{equation}

To proceed, we apply Proposition \ref{prop5.3}. First recall that, according
to Karamata's theorem, any positive slowly varying function $h(x)$
defined in $x \geq0$ has a representation
\[
h(x) = c(x)\exp\biggl\{\int_{x_0}^x
\frac{w(y)}{y}\mrmd y \biggr\},\qquad x \geq x_0,
\]
where $x_0>0$, $c(x) \rightarrow1$, and $w(x) \rightarrow0$,
as $x \rightarrow\infty$. For $x_0 = \min_{n \geq1} b_n$,
$1 \leq|t| \leq\ep b_n$, where $0 < \ep\leq1$ is fixed, this
representation implies that with some constant $c_0 > 0$
\[
\frac{h(b_n/|t|)}{h(b_n)} \geq c_0 |t|^{-\gamma} \qquad\mbox{with }
\gamma= \gamma(\ep) = \sup_{y \geq1/\ep} \bigl|w(y)\bigr|.
\]
Hence, from (\ref{equ5.7})--(\ref{equ5.8})
\[
\bigl|f_1(t/b_n)\bigr| = \exp\bigl\{-c|t|^{\alpha}
b_n^{-\alpha} h\bigl(b_n/|t|\bigr)\bigr\} \leq\exp\bigl
\{-c_1 |t|^{\alpha- \gamma}/n\bigr\}
\]
with some constant $c_1 > 0$.

We choose $\ep>0$ to be small enough so that $\gamma< \alpha- \delta$.
Now, applying the above estimate in (\ref{equ7.5}), we get in the region
$1 \leq|t| \leq\ep b_n$
\begin{eqnarray*}
\bigl|\tilde{f}_n(t)\bigr| & \leq& \frac{1}{1-\delta_n} \bigl(\bigl|f_1(t/b_n)\bigr|
+ \delta_n \bigr)
\\
& \leq& \frac{1}{1-\delta_n} \bigl(\exp\bigl\{-c_1
|t|^{\alpha- \gamma}/n\bigr\} + \delta_n \bigr).
\end{eqnarray*}
One can simplify the right-hand side by noting that
$\frac{c_1 |t|^{\alpha- \gamma}}{n} \leq
\frac{c_1 b_n^{\alpha- \gamma}}{n} < K$
with some constant~$K$. Using $\log x \leq x-1$ ($x > 0$) and
$\RMe^{-x} \leq1 - \frac{1}{K} (1 - \RMe^{-K}) x$, for $0 \leq x \leq K$,
we then have
\begin{eqnarray*}
\log\bigl(\exp\bigl\{-c_1 |t|^{\alpha- \gamma}/n\bigr\} +
\delta_n \bigr) & \leq& \exp\bigl\{-c_1
|t|^{\alpha- \gamma}/n\bigr\} + \delta_n - 1
\\
& \leq& -\frac{1 - \RMe^{-K}}{K} \frac{c_1 |t|^{\alpha- \gamma}}{n} +
\delta_n \\
&\leq&
\frac{c_2}{n} - \frac{c_3 |t|^{\alpha- \gamma}}{n}
\end{eqnarray*}
with positive constants $c_j$. As a result,
\[
\bigl|\tilde{f}_n(t)\bigr| \leq\exp\biggl\{\frac{1}{n}
\bigl(c_4 - c_5 |t|^{\alpha- \gamma} \bigr) \biggr\}
\]
with some other positive constants $c_4$ and $c_5$ (independent of $n$).
It remains to raise this inequality to the power $k$, and (\ref
{equ7.4}) follows.
\end{pf}

We will now develop a few applications of Lemmas \ref{lem7.1} and \ref
{lem7.2} using the
assumption \mbox{$I(p) < \infty$}. The latter forces $p$ to have bounded
variation and vanish at infinity. Hence,
%
%
\begin{equation}\label{equ7.6}
\|r_n\|_{\mathrm{TV}} = \|\tilde{p}_n\|_{\mathrm{TV}}
= b_n(1 - \delta_n)^{-1} \|p 1_{\{|x| \leq b_n\}}
\|_{\mathrm{TV}} \leq b_n (1 - \delta_n)^{-1}
\sqrt{I(p)}.
\end{equation}
Using the inequality (\ref{equ3.2}), we see that the characteristic functions
of $\tilde{p}_n$ and of the centered density
$r_n(x) = \tilde{p}_n(x+d_n)$ satisfy
%
%
\begin{equation}\label{equ7.7}
\bigl|\psi_n(t)\bigr| = \bigl|\tilde f_n(t)\bigr| \leq\frac{cb_n}{|t|}
\qquad(t \neq0)
\end{equation}
with some constant $c = c(p)$, depending on $p$, only.

%
\begin{corollary}\label{corol7.3}
If $I(p) < \infty$, then under the assumptions
of Lemma \ref{lem7.2} with $k \geq4$, we have with some constant $C$ depending
on $p,\delta,\eta$, only,
%
%
\begin{eqnarray}\label{equ7.8}
\int_{-\infty}^\infty\bigl(1+|t|\bigr) \bigl|\psi_n^k(t)\bigr|
\mrmd t &\le& C,
\\
\label{equ7.9}
\int_{-\infty}^\infty t^2 \bigl|\bigl(
\psi_n^k\bigr)'(t)\bigr|^2\mrmd t &\le& C.
\end{eqnarray}
\end{corollary}

\begin{pf} We have $(\psi_n^k)'(t) = k\psi_n'(t) \psi_n(t)^{k-1}$,
while by (\ref{equ7.2}),
\[
\int_{-\infty}^\infty t^2 \bigl|
\psi_n'(t)\bigr|^2 \bigl|\psi_n(t)\bigr|^{2(k-1)}
\mrmd t \le\frac{C^2}{n^2} \int_{-\infty}^\infty
t^4 \bigl|\psi_n(t)\bigr|^{2(k-1)}\mrmd t.
\]
To estimate the last integral, first we use (\ref{equ7.4}) which gives
\[
\int_{|t| \le\ep b_n} t^4 \bigl|\psi_n(t)\bigr|^{2(k-1)}
\mrmd t \le C.
\]
For the complementary region $|t| > \ep b_n$, note that
\[
\tilde f_n(b_n t) = \frac{1}{1-\delta_n} \int
_{-b_n}^{b_n} \RMe^{\RMi tx} p(x)\mrmd x,
\]
which shows that these functions are separated from 1 uniformly in $n$
in $|t| \geq\ep$. (This can easily be seen by using general separation
bounds for characteristic functions which are discussed in~\cite{B-C-G2}.)
Thus,
\[
\sup_{|t|\ge\varepsilon}\bigl|\psi_n(b_nt)\bigr| = \sup
_{|t|\ge\varepsilon
}\bigl|\tilde{f}_n(b_nt)\bigr| \leq
\RMe^{-c}
\]
for some constant $c>0$ independent of $n$. In addition, by (\ref{equ7.7}),
\[
t^4 \bigl|\psi_n(b_nt)\bigr|^6 \le
\frac{c}{t^2}
\]
with some other constant. Hence,
\[
\int_{|t|\ge\varepsilon b_n} t^4 \bigl|\psi_n(t)\bigr|^{2(k-1)}
\mrmd t \le b_n^5 \RMe^{-2c(k-4)}\int_{|t|\ge\varepsilon}
t^4 \bigl|\psi_n(b_nt)\bigr|^6\mrmd t \le
Cb_n^5 \RMe^{-2ck}.
\]
The last expression is exponentially small with respect to $n$ by the
constraint on $k$, and we arrive at (\ref{equ7.9}). The first
inequality (\ref{equ7.8}),
which is simpler, is proved similarly.
\end{pf}

\section{Boundedness of Fisher information. Proof of Theorem \texorpdfstring{\protect\ref{theo1.1}}{1.1}}\label{sec8}

In this section, we complete the last step in the proof of Theorem \ref
{theo1.1}.
Keeping the same notations as in the previous sections and recalling
Lemma \ref{lem6.1}, we only need the following lemma.

%
\begin{lemma}\label{lem8.1}
Assume that $Z_n \Rightarrow Z$ weakly in
distribution, where $Z$ has a non-extremal stable law. If
$I(Z_{n_0}) < \infty$ for some $n_0$, then
$\sup_{n \geq n_0} I(Z_n) < \infty$.
\end{lemma}

In the normal case, when $X_1$ has a finite second moment, the assertion
immediately follows from Stam's inequality (\ref{equ2.3}). In view of
Lemma \ref{lem6.1},
we therefore obtain Barron--Johnson theorem, that is,
$I(Z_n\|Z) \rightarrow0$.
Thus, we may focus on the case $0 < \alpha< 2$.

To simplify the argument and the notations, we assume $n_0 = 1$
(otherwise, mild modifications connected with the binomial decomposition
are only needed). Thus, let $I(p) < \infty$, where $p$ is the density
of $X_1$. As in the previous section, we denote by $p_n$ the density
of $\tilde Z_n = Z_n + a_n$ and assume that $Z_n \Rightarrow Z$ weakly
in distribution, where $Z$ has a non-extremal stable law.

By Stam's inequality (\ref{equ2.3}),
\[
I(Z_n) \leq\frac{b_n^2}{n} I(p).
\]
Although the right-hand side tends to infinity, as $n \rightarrow
\infty$,
this inequality may be used for small values of $n$, and here it will be
sufficient to show that $\sup_{n \geq n_0} I(Z_n) < \infty$ for some
$n_0$.

Our basic tool is the binomial decomposition (\ref{equ7.1}) of the previous
section. Note that, by the convexity of the $I$-functional,
%
%
\begin{equation}\label{equ8.1}
I(p_n) \le\sum_{k=0}^n \pmatrix{n
\cr k} (1-\delta_n)^k\delta_n^{n-k}
I \bigl(\tilde{p}_n^{k*}*\tilde{q}_n^{(n-k)*}
\bigr),
\end{equation}
so it will be sufficient to properly estimate the terms in this sum.
To this aim, we fix a number $\eta\in(0,1)$ and distinguish two cases.

%
\begin{lemma}\label{lem8.2}
If $k \leq n-3$, then
%
%
\begin{equation}\label{equ8.2}
I \bigl(\tilde{p}_n^{k*}*\tilde{q}_n^{(n-k)*}
\bigr) \le C(nb_n)^2 I(p)
\end{equation}
with some constant $C$ depending on $p$, only.
\end{lemma}

\begin{pf}
By the monotonicity property (\ref{equ2.2}),
$I(\tilde{p}_n^{k*}*\tilde{q}_n^{(n-k)*})\le I(\tilde{q}_n^{(n-k)*})$.
On the other hand, by Proposition \ref{prop3.2}, if $n-k \geq3$,
\[
I\bigl(\tilde{q}_n^{(n-k)*}\bigr) \le\tfrac{1}{2} \bigl(
\bigl\|\tilde{q}_n^{[(n-k)/3]*}\bigr\|_{\mathrm{TV}}^2 + 2 \bigl\|
\tilde{q}_n^{[(n-k)/3]*}\bigr\|_{\mathrm{TV}} \cdot\bigl\|
\tilde{q}_n^{n-k-2[(n-k)/3]*}\bigr\|_{\mathrm{TV}} \bigr).
\]
But the total variation norm decreases when taking convolutions,
so that $\|\tilde{q}_n^{s*}\|_{\mathrm{TV}}\le\|\tilde{q}_n\|_{\mathrm{TV}}$
($s = 1,2,\ldots$). Hence,
\[
I\bigl(\tilde{q}_n^{(n-k)*}\bigr) \le\tfrac{3}{2} \|
\tilde{q}_n\|_{\mathrm{TV}}^2.
\]
In turn, by means of the inequality $\|p\|_{\mathrm{TV}} \leq\sqrt{I(p)}$
(Proposition \ref{prop3.1}), we have
\[
\|\tilde{q}_n\|_{\mathrm{TV}} = b_n
\delta_n^{-1} \|p 1_{\{|x|>b_n\}}\|_{\mathrm{TV}} \leq
b_n\delta_n^{-1} \|p\|_{\mathrm{TV}} \leq
b_n\delta_n^{-1}\sqrt{I(p)},
\]
where we used the property $p(-\infty)=p(\infty)=0$ for the first
inequality. Thus
\[
I \bigl(\tilde{p}_n^{k*}*\tilde{q}_n^{(n-k)*}
\bigr) \le\tfrac
{3}{2} \bigl(\sqrt{I(p)} b_n
\delta_n^{-1} \bigr)^2.
\]
Recalling that $\delta_n \sim\frac{c}{n}$, Lemma \ref{lem8.2} is proved.
\end{pf}

%
\begin{lemma}\label{lem8.3}
If $15 \leq\eta n \le k \leq n$, then
%
%
\begin{equation}\label{equ8.3}
I \bigl(\tilde{p}_n^{k*}*\tilde{q}_n^{(n-k)*}
\bigr) \le C
\end{equation}
with some constant $C$ depending on $p$ and $\eta$, only.
\end{lemma}

\begin{pf} Again appealing to the monotonicity of the Fisher
information, we will use the bound
\[
I\bigl(\tilde{p}_n^{k*}*\tilde{q}_n^{(n-k)*}
\bigr)\le I\bigl(\tilde{p}_n^{k*}\bigr).
\]
Thus, involving the centered density $r_n(x) = \tilde{p}_n(x+d_n)$
with the characteristic function $\psi_n$ (as in the previous section),
it suffices to show that
%
%
\begin{equation}\label{equ8.4}
I\bigl(r_n^{k*}\bigr) = I\bigl(\tilde{p}_n^{k*}
\bigr) \le C.
\end{equation}

Assume first that $\eta_0 n \leq k \le n$, where $0 < \eta_0 < \eta$.
Since $\|r_n\|_{\mathrm{TV}} \leq Cb_n \sqrt{I(p)} < \infty$
(see (\ref{equ7.6}) and Proposition \ref{prop3.2}), the convolution
powers $r_n^{k*}$
have finite Fisher information, whenever $k \geq3$. In view of the
bound (\ref{equ7.7}) on the characteristic functions, we may invoke inversion
formulas like in (\ref{equ5.5})--(\ref{equ5.7}) to write, for any $x \in
\R$,
%
%
\begin{eqnarray}\label{equ8.5}
r_n^{k*}(x) & = & \frac{1}{2\uppi}\int
_{-\infty}^{\infty} \RMe^{-\RMi tx}\psi_n(t)^{k}
\mrmd t,
\\
\label{equ8.6}
\bigl(r_n^{k*}\bigr)'(x) & = &
\frac{1}{2\uppi}\int_{-\infty}^{\infty} \RMe^{-\RMi tx}
(-\RMi t)\psi_n(t)^{k}\mrmd t,
\\
\label{equ8.7}
r_n^{k*}(x) + x\bigl(r_n^{k*}
\bigr)'(x) & = & - \frac{1}{2\uppi} \int_{-\infty}^{\infty}
\RMe^{-\RMi tx} tk\psi_n(t)^{k-1}\psi_n'(t)
\mrmd t,
\end{eqnarray}
where for reasons of integrability it is safer to assume that $k \geq5$.

Corollary \ref{corol7.3} tells us that the Fourier transforms in (\ref
{equ8.5}) and (\ref{equ8.7})
are well defined for square integrable functions whose $L^2$-norms are
bounded by a constant independent of $k$ and $n$. Hence, the same is
true for
\[
x \bigl(r_n^{k*}\bigr)'(x) = -
\frac{1}{2\uppi} \int_{-\infty}^{\infty} \RMe^{-\RMi tx}
\bigl(\psi_n(t)^{k} + tk\psi_n(t)^{k-1}
\psi_n'(t) \bigr)\mrmd t,
\]
and we may write
%
%
\begin{equation}\label{equ8.8}
\bigl|\bigl(r_n^{k*}\bigr)'(x)\bigr| \leq
\frac{u_{nk}(x)}{|x|}
\end{equation}
with
%
%
\begin{equation}\label{equ8.9}
\|u_{nk}\|_2^2 = \int_{-\infty}^\infty
u_{nk}(x)^2\mrmd t \leq C.
\end{equation}
Moreover, according to (\ref{equ7.8}), $L^1$-norms of the functions
$(-\RMi t)\psi_n(t)^{k}$ in (\ref{equ8.6}) are also bounded
by a constant independent of $k$ and $n$. Hence,
\[
\sup_x \bigl|\bigl(r_n^{k*}
\bigr)'(x)\bigr|\le C
\]
for all $n$ and $\eta_0 n \leq k\le n$. As a result, (\ref{equ8.8}) may be
sharpened to
\[
\bigl|\bigl(r_n^{k*}\bigr)'(x)\bigr| \leq
\frac{u_{nk}(x)}{1+|x|}
\]
with some functions $u_{nk}$ satisfying (\ref{equ8.9}). By applying Cauchy's
inequality, the latter immediately implies that
%
%
\begin{equation}\label{equ8.10}
\bigl\|r_n^{k*}\bigr\|_{\mathrm{TV}} = \int_{-\infty}^\infty\bigl|
\bigl(r_n^{k*}\bigr)'(x)\bigr|\mrmd x \le
C' \|u_{nk} \|_2 \le C,
\end{equation}
where the resulting constant $C$ may depend on $p$ and $\eta_0$
(by choosing, e.g., $\delta= \alpha/2$ in the previous auxiliary
lemmas of the previous section).

We now apply Proposition \ref{prop3.2} to convolutions of any three densities
$r_n^{k*}$, as above. That is, if $\eta_0 n \leq k_j \le n$ and
$k_j \geq5$ ($j=1,2,3$), we obtain by (\ref{equ3.4}) and (\ref
{equ8.10}) that
%
%
\begin{equation}\label{equ8.11}
I\bigl(r_n^{(k_1 + k_2 + k_3)*}\bigr) \leq\tfrac{3}{2}
C^2.
\end{equation}

Starting with $k \geq15$, put $k_1 = k_2 = [\frac{k}{3}]$,
$k_3 = k - (k_1+k_2)$, so that $k_j \geq5$. Also, if $k \geq\eta n$,
we have $k_j \geq[\frac{\eta n}{3}] \geq\frac{\eta n}{6}$.
Hence, we may choose $\eta_0 = \frac{\eta}{6}$, and thus (\ref{equ8.11})
implies (\ref{equ8.3})--(\ref{equ8.4}).
\end{pf}

\begin{pf*}{Proof of Lemma \ref{lem8.1}} In the case $15 \leq\eta n \le
n-3$, we may
combine Lemmas \ref{lem8.2} and \ref{lem8.3} to get from (\ref{equ8.1})
the following.
With some constant
$C = C(p,\eta)$, depending on $\eta$ and the density $p$ via $I(p)$
and the constant $c$ in $\delta_n \sim\frac{c}{n}$,
\begin{eqnarray*}
I(p_n) & \le& C(nb_n)^2 \sum
_{0 \leq k < \eta n} \pmatrix{n\cr k} (1-\delta_n)^k
\delta_n^{n-k} + C\sum_{\eta n \leq k \leq n} \pmatrix{n
\cr k}(1-\delta_n)^k\delta_n^{n-k}
\\
& \le& C(nb_n)^2 \cdot2^n
\delta_n^{(1-\eta) n} + C \le C',
\end{eqnarray*}
where the last inequality holds for all sufficiently large $n$
(by using $\delta_n \sim\frac{c}{n}$) with, for example,
$\eta= \frac{1}{2}$. Lemma \ref{lem8.1} and therefore Theorem \ref
{theo1.1} are now proved.
\end{pf*}

%
\begin{remark}\label{rem8.4}
Finally, let us comment on the conditions \textup{(a)--(b)}
from the \hyperref[sec1]{Introduction}. In view of the general bound
(\ref{equ3.3}),
\textup{(a)} is always necessary for the finiteness of $I(Z_n)$ with some $n$.
Since \textup{(b)} is weaker than \textup{(a)}, we need explain the opposite direction.

If $1 < \alpha\leq2$, then $X_1$ has finite first absolute moment
$C = \E|X_1|$. Hence, under (\ref{equ1.6}), the condition (\ref
{equ3.5}) is fulfilled
and thus the bound (\ref{equ3.6}) is applicable to all $Z_n$ with
$n \geq(\nu+ 2)/2$. More precisely, denoting by $g_n(t) = f_1(t)^n$
the characteristic function of $S_n = X_1 + \cdots+ X_n$, we have
\[
\bigl|\bigl(tg_n(t)\bigr)'\bigr| \leq\bigl|g_n(t)\bigr| + |t|
\bigl|g_n'(t)\bigr| \leq\bigl(1 + Cn|t|\bigr) \bigl|f_1(t)\bigr|^{n-1},
\]
thus $S_n$ has a density $\rho_n(x)$ whose total variation norm satisfies
\[
\|\rho_n\|_{\mathrm{TV}}^4 \leq\int
_{-\infty}^\infty t^2 \bigl|f_1(t)\bigr|^{2n}
\mrmd t \int_{-\infty}^\infty\bigl(1 + Cn|t|\bigr)^2
\bigl|f_1(t)\bigr|^{2(n-1)} \mrmd t < \infty.
\]
By Proposition \ref{prop3.2}, we get $I(S_{3n}) < \infty$.
\end{remark}

\section*{Acknowledgements}

Research partially supported by NSF Grant DMS-11-06530,
Simons Foundation and SFB 701.
We thank the referees for careful reading of
the manuscript and valuable comments.


%

\printhistory

\end{document}